%% file: main.tex
\documentclass{article}

\usepackage[utf8]{inputenc}
\usepackage[T1]{fontenc}
\usepackage{amsmath,amssymb,amsthm,bbm,float,graphicx,geometry,lmodern,mathtools,parskip,setspace}
\usepackage[colorlinks, linkcolor = blue!80!black, citecolor = blue!80!black,breaklinks, pdfauthor={Benedikt Jahnel, Daniel Kamecke}]{hyperref}
\usepackage{mathrsfs}
\usepackage{dsfont}
\usepackage{enumerate}
\usepackage{nicefrac}
\usepackage{todonotes}
\usepackage{comment}
\usepackage{titling}
\usepackage{fullpage}
\usepackage{comment}
\usepackage{orcidlink}
\usepackage{enumitem}
\usepackage{graphicx}
\usepackage{placeins}
\usepackage{float}
\usepackage{subfig}

\usepackage{tikz}
\usetikzlibrary{backgrounds}
\usetikzlibrary{patterns}
\usetikzlibrary{positioning, shapes.geometric}

\usepackage{caption}
\usepackage{graphicx}
\graphicspath{{Images/}}
\allowdisplaybreaks

\newcommand{\N}{\mathbb{N}}

\newtheorem{theorem}{Theorem}[section]

\newtheorem{corollary}[theorem]{Corollary}
\newtheorem{prop}[theorem]{Proposition}

\theoremstyle{definition}

\newtheorem{example}[theorem]{Example}
\newtheorem{remark}[theorem]{Remark}

\usepackage{nameref}

\makeatletter
\let\orgdescriptionlabel\descriptionlabel
\renewcommand*{\descriptionlabel}[1]{%
  \let\orglabel\label
  \let\label\@gobble
  \phantomsection
  \edef\@currentlabel{#1}%
  \let\label\orglabel
  \orgdescriptionlabel{#1}%
}

\newcommand{\R}{\mathbb{R}}
\renewcommand{\N}{\mathbb{N}}
\newcommand{\Q}{\mathbb{Q}}
\newcommand{\E}{\mathbb{E}}

\newcommand{\M}{\mathbb{M}}
\renewcommand{\P}{\mathbb{P}}

\newcommand{\supp}{\operatorname{supp}}

\newcommand{\compl}{\mathsf{c}}
\newcommand*\diff{\mathop{}\!\mathrm{d}}

\pretitle{\centering\LARGE\scshape}
 \posttitle{\vskip 0.75cm}

 \predate{\vskip 0.75 cm \centering\large}
 \postdate{\par}

\usepackage[style = numeric, sorting=nyt, url = true, abbreviate=false, maxbibnames=9, sortcites=true, doi = false, backend = biber, giveninits = true, isbn=false]{biblatex}
\renewbibmacro{in:}{\ifentrytype{article}{}{\printtext{\bibstring{in}\intitlepunct}}}
\bibliography{refs/refs.bib}

\title{Phase transitions for the Widom--Rowlinson model in random environments}

\thanksmarkseries{arabic}

\author{
Benedikt Jahnel \orcidlink{0000-0002-4212-0065}\thanks{Technische Universit\"at Braunschweig, Universit\"atsplatz 2, 38106 Braunschweig, Germany}\thanksgap{0.4ex} \thanks{Weierstrass Institute for Applied Analysis and Stochastics, Mohrenstr.\ 39, 10117 Berlin, Germany} \\ benedikt.jahnel@tu-braunschweig.de
\and
Daniel Kamecke \orcidlink{}\thanksmark{1}\\daniel.kamecke@tu-braunschweig.de
}

\date{May 8, 2025}

\begin{document}

\maketitle
\begin{spacing}{0.9}
\begin{abstract}  
We establish non-uniqueness regimes for the infinite-volume two-colored Widom--Rowlinson model based on inhomogeneous Poisson point processes with locally finite intensity measures featuring percolation. As an application, we provide almost-sure phase-transition results for the Widom--Rowlinson model based on translation-invariant and ergodic Cox point processes with stabilizing and non-stabilizing directing measures.  

\smallskip
\noindent\footnotesize{{\textbf{AMS-MSC 2020}: 60K35, 05C80}

}
\end{abstract}
\end{spacing}


\section{Introduction} \label{sec:Intro}
The continuum Widom--Rowlinson model is arguably the most paradigmatic Gibbs point process for which a phase transition of non-uniqueness of solutions for the associated DLR equations has been established, see~\cite{Ru71,ChChKo95,GiLeMa95, GeHa95,DeHo21,DeHo15}. Designed as a toy model for liquid-vapor transitions in the continuum, see~\cite{WiRo70}, it has enjoyed much attention in the last decades, in part, due to its accessibility to tools in stochastic geometry. Indeed, using the random-cluster representation for the symmetric two-colored model, the phase transition can be established in regimes where the underlying Poisson point process features sufficiently strong percolation. Here, roughly speaking, percolation guarantees that information from an infinitely far-away boundary condition can be felt in the particle density. 

The purpose of this short note is to highlight that this general strategy remains intact also in cases where the intensity of the underlying Poisson point process is not necessarily stationary with respect to shifts in $\R^d$. We only require that the intensity measure is locally finite and that the associated Poisson--Boolean model is in the supercritical percolation regime for a slightly increased intensity. 

Our interpretation is that the intensity measure represents an environment that may feature inhomogeneities, like a material with impurities, and then the questions becomes, which environments do feature non-trivial supercritical percolation phases. Apart from simple examples, e.g., where the inhomogeneous Poisson point process can be coupled with a supercritical homogeneous one, we apply our result to a class of models where the environment is random. 

More precisely, we consider Cox point processes, that is Poisson point processes in random environment, where the random environment measure, i.e., intensity measure, has a stationary distribution. In this case, the Cox point process is again stationary and continuum-percolation results have been established in recent years~\cite{HiJaCa19,JaToCa22,HiJaMu22} under certain mixing conditions, but also in strongly correlated settings~\cite{JaJhVu23}. With the help of these results, we are then in the position to establish almost-sure phase transition results for the Widom--Rowlinson model in large classes of random environments. 

Let us note that, for simplicity, we focus here only on the classical Widom--Rowlinson model with two possible colors and balls with fixed radius and hard-core interaction. Many generalizations of this model can be considered in future work, such as Potts models with additional particle interaction~\cite{GeHa95}, lattice versions~\cite{HiTa04,Ku19}, additional marks that may represent random radii~\cite{DeHo19}, and many more. 

Let us also mention that the Widom--Rowlinson model in random environment can be viewed as a continuum version of the random-field Ising model~\cite{FyMaPiSo18,BrKu88,Ch24,DiLiXi24} in the sense that the role of the iid perturbations in the field term in the Hamiltonian of the random-field Ising model is played by the random intensity measure. In particular, our setting then also covers situations in which the random field is correlated. 

The manuscript is organized as follows. In the next Section~\ref{sec:Setting_Main_Result}, we present our setting and the main results, first for fixed environments and then for random environments. In Section~\ref{sec:Examples} we present a number of examples that highlight the applicability of our results. Finally, in Section~\ref{sec:Proofs} we give proofs for the main results. 

\section{Setting and main result}\label{sec:Setting_Main_Result}
Before stating our main result in Section~\ref{sec:Main_Results}, we briefly introduce some notation in Section~\ref{sec:Setting}. 
\subsection{Two-colored Widom--Rowlinson models in inhomogeneous environments}\label{sec:Setting}
Let $\mathbb{M}=\R^d\times \{-,+\}$ be the \textit{marked space}, where the first variable of $(x,c)\in \mathbb{M}$ represents the spatial {\em position}  and $c$ its {\em mark}. We also call the marks $+$ and $-$ \textit{colors}. Define, for $\Lambda\subset \mathbb{M}$, the \textit{space of marked point configurations} supported on $\Lambda$ by  
\begin{align}\label{Def_Gamma}
    \Gamma_{\Lambda}=\{\pmb{\omega}\subset\Lambda\,|\, \forall \Delta\Subset \Lambda\colon    |\pmb{\omega}_{\Delta}|<\infty \},
\end{align}
where $\Delta\Subset \Lambda$ is a bounded subset, $\pmb{\omega}_\Delta=\pmb{\omega}\cap \Delta$ and $|\pmb{\omega}|$ is the cardinality of $\pmb{\omega}$. The definition of $\Gamma_\Lambda$ in \eqref{Def_Gamma} can also be used for $\Lambda\subset \R^d$. Furthermore, we denote by $\mathfrak{F}_\Lambda$ the usual Borel $\sigma$-algebra on $\Gamma_\Lambda$, see, e.g., \cite[Chapter~9.1]{DaVe07}. In the following, we sometimes call $\sigma$ an \textit{environment}, when it is a locally finite measure on $\R^d$. Given $\sigma$, we are also interested in the corresponding measure $\bar \sigma $ on $\M$ which is symmetric in both colors, i.e., $\bar\sigma=\sigma\otimes (\delta_{-}+\delta_+)$, where $\delta_c$ is the Dirac measure and $\otimes$ denotes the product measure. For an environment $\sigma$ and $\Lambda\subset \R^d$, let $\Pi_{\sigma,\Lambda}$ be the law of a Poisson point process on $\Gamma_{\R^d}$ with intensity measure $\sigma_\Lambda(A)=\sigma(A\cap \Lambda)$. Analogously, define $\Pi_{\bar\sigma,\Lambda}$ for $\Lambda\subset \M$ as a the law of a Poisson point process on $\Gamma_\M$. For more information on Poisson point processes, see, e.g., \cite{GuPe17}. 

Next, fix $a>0$ and define the \textit{viable configurations} as
\begin{align*}
    \pmb{V}=\{\pmb{\omega}\in\Gamma_{\mathbb{M}}\,|\, B_a(\pmb{\omega}^+)\cap B_a(\pmb{\omega}^-)=\emptyset \},
\end{align*}
where we used the notation $\pmb{\omega}^{\pm}=\{x\in \R^d\,|\, (x,\pm)\in \pmb{\omega}\}$ for $\pmb{\omega}\subset \M$ as well as $B_a(\omega)=\bigcup_{x\in \omega} B_a(x)$ for $\omega \subset \R^d$, where $B_a(x)$ is the radius-$a$ ball centered at $x$. Consider, for $\pmb{\omega} \subset \mathbb{M}$, the \textit{specification of the Widom--Rowlinson model} given by 
\begin{align*}
    \gamma^\sigma_\Lambda(\diff \pmb{\zeta} \,|\, \pmb{\omega})=\frac{1}{Z^\sigma_\Lambda(\pmb{\omega})}\mathds{1}_{\pmb{V}}(\pmb{\zeta}\cup\pmb{\omega}_{\Lambda ^\compl}) \, \Pi_{ \bar{\sigma},\Lambda }(\diff \pmb{\zeta}) \quad \text{for }\Lambda\Subset \M  ,
\end{align*}
where $Z^\sigma_\Lambda(\pmb{\omega})$ normalizes the expression. Note that $0<Z^\sigma_\Lambda(\pmb{\omega})\le 1$, for $\Lambda\Subset \M$, since $\sigma$ is assumed to be locally finite. Further, note that in this definition, we do not use $\pmb{\omega}\in \Gamma_{\mathbb{M}}$ but allow all $\pmb{\omega}\subset \mathbb{M}$.
However, when restricting to $\pmb{\omega} \in \Gamma_{\mathbb{M}}$, the specification of the Widom--Rowlinson model $\gamma^\sigma_\Lambda$ is a probability kernel from $\Gamma_{\mathbb{M}}$ to $\Gamma_{\mathbb{M}}$. We say that $\mu$ is a \textit{Gibbs measure of the Widom--Rowlinson model in the environment~$\sigma$ with parameter $a$}, or $\mu\in \mathcal{G}_{\sigma,a}$, if, for all $\Lambda\Subset\M $ and $A\in \mathfrak{F}_\M$,
\begin{align*}
\mu(A)=\int\gamma_\Lambda^\sigma(A\,|\,\pmb{\zeta})\,\mu(\diff \pmb{\zeta}).
\end{align*}
A central theme within the study of Gibbs measures is the question of non-uniqueness, i.e., $|\mathcal{G}_{\sigma,a}|>1$. This is due to the idea that at those $\sigma$ and $a$ the system could, in principle, change from one state to another, which can be interpreted as a phase transition. For more information on Gibbs measures, see \cite{Ge11,De19,Ja18}. After having established the necessary setup, we now present our results. 
\subsection{Main results}\label{sec:Main_Results}
The proof of the celebrated classical phase-transition result for the Widom--Rowlinson model connects the non-uniqueness of Gibbs measures with percolation. The notion of percolation is not unified in the literature. For our purposes it is practicle to say that $\omega \in \Gamma_{\R^d}$ is $a$\textit{-percolating} if there is an infinite connected component in $B_a(\omega)$. Although the connection between percolation and non-uniqueness of the Widom--Rowlinson model is usually stated only for the case where the environment is given by multiples of the Lebesgue measure~\cite{ChChKo95,DeHo19}, it holds true also for more general environments. 
\begin{theorem}
    \label{prop:Non-Random_Environment}
    Let $a>0$ and $\sigma$ a locally finite measure on $\R^d$. There is a $\tau>0$, only dependent on the dimension $d\ge 1$, such that
\begin{align}\label{Cond:Percolation}
        \Pi_{\tau \sigma, \R^d}(\{\omega \,|\,  \omega \text{ is }a\text{-percolating}\})>0\quad \implies \quad |\mathcal{G}_{\sigma,a}|>1.
    \end{align}
\end{theorem}
The proof of Theorem~\ref{prop:Non-Random_Environment} is postponed to Section~\ref{sec:Proofs}. The following corollary, first presented in~\cite{ChChKo95,GiLeMa95}, follows by percolation results for the Poisson point process, see, e.g., \cite[Page~52]{Me96}. 
\begin{corollary}\label{Cor:Lebesgue_Phase_Transition}
    Let $d\ge2$, $a>0$ and $\lambda$ the Lebesgue measure on $\R^d$. Then, for all $z>0$ large enough, $|\mathcal{G}_{z\lambda,a}|>1$.
\end{corollary}
Here, we focus on sufficient conditions for non-uniqueness. Conversely, one might be interested in a criterion for uniqueness. In the case $\sigma=\lambda$, it is proved that for $z$ small enough it holds $|\mathcal{G}_{z\lambda, a}|=1$, see, e.g., \cite[Theorem~10.2]{GeHaMa01}. 
We comment on the implication that absence of percolation should guarantee uniqueness, in our more general setting of Theorem~\ref{prop:Non-Random_Environment}, below in  Remark~\ref{Rem:Uniqueness}.

Besides multiples of the Lebesgue measure, a rich class of environments that feature percolation, and hence non-uniqueness of Gibbs measures, is obtained when replacing $\sigma$ by a random environment $\Sigma$. More precisely, $\Sigma$ is assumed to be a \textit{random measure}, i.e., a random variable $\Sigma\colon(\Omega,\mathcal{F}, \P)\rightarrow(\mathcal{M},\mathscr{F})$, where $(\Omega,\mathcal{F}, \P)$ is a probability space and $\mathcal{M}$ is the space of Borel measures on $\R^d$ with Borel $\sigma$-algebra $\mathscr{F}$ as in \cite[Chapter~9.1]{DaVe07}. The distribution of a \emph{Cox point process w.r.t.~$z\Sigma$} is given by $\Q_z(A)= \E_\P\left[\Pi_{z\Sigma, \R^d}(A)\right]$. Then, non-uniqueness of the Widom--Rowlinson model in the environment $\Sigma$ is guaranteed by percolation of a corresponding Cox point processes.
\begin{corollary}\label{thm:Random_Environment}
    Let $a>0$ and $\Sigma$ a random measure. There is a $\tau>0$, only dependent on the dimension $d\ge 1$, such that
\begin{align}\label{Cond:Percolation2}
        \Q_\tau (\{\omega \,|\,  \omega \text{ is }a\text{-percolating}\})=1\quad &\implies \quad \P( |\mathcal{G}_{\Sigma,a}|>1)=1,\\
        \Q_\tau (\{\omega \,|\,  \omega \text{ is }a\text{-percolating}\})>0\quad &\implies \quad \P( |\mathcal{G}_{\Sigma,a}|>1)>0.
    \end{align}
\end{corollary}
As percolation of Cox point processes is studied in literature already, see \cite{HiJaCa19,JaJhVu23}, we readily get non-uniqueness results for the Widom--Rowlinson model in random environment. In the following, we showcase some examples of Cox point processes, where the existence of percolation is proved.

\section{Examples}\label{sec:Examples}
This section is devoted to examples of environments such that the Widom--Rowlinson model exhibits more than one Gibbs measure. Already the class of non-random environments as considered in Theorem~\ref{prop:Non-Random_Environment} is quite rich. 

\begin{figure}[!htpb]
\centering
\includegraphics[width=0.9\columnwidth]{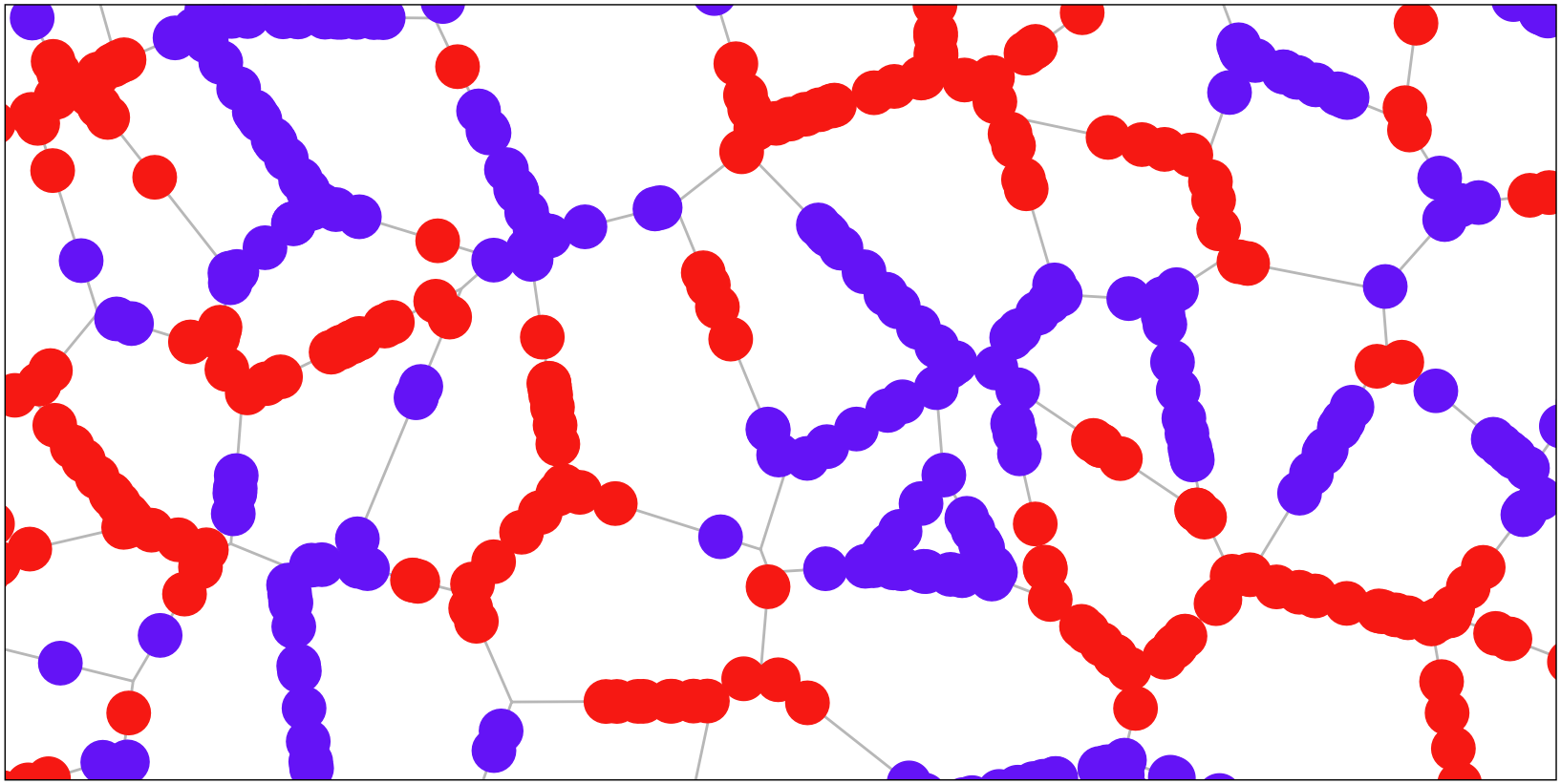}
\caption{Realization of the two-colored Widom--Rowlinson model based on Cox point process with intensity measure given by a Poisson--Voronoi tessellation.}
\label{Pix_3}
\end{figure}

\begin{example}[Non-random measures]
    A measure $\sigma$ such that $\sigma(A)\ge \alpha \lambda(A)$, where $\alpha$ large and $\lambda$ the Lebesgue measure satisfies the assumptions in Theorem~\ref{prop:Non-Random_Environment}. In fact, the domination $\sigma(A)\ge \alpha \lambda(A)$ can be relaxed to the case where $\sigma(B)=0$ is allowed for a lot of $B\subset \R^d$ since most arguments used to prove continuum percolation are rather flexible towards adaptations of the Lebesgue measure~\cite{Me96,BoRi06,JaKo20}. 
\end{example}
Another rich class of suitable environments is provided by considering random environments as used in Corollary~\ref{thm:Random_Environment}. 
\begin{figure}[!htpb]
    \centering
    \subfloat[\label{Pix_1}Realization of Example~\ref{Example:Absolutely_Continuous}~\ref{Example:Random_closed_set}, where $\lambda_1=1$, $\lambda_2=0$ and $\Xi$ is given by an independent Poisson--Boolean model, i.e., $\Xi=\bigcup_{i\in I}B_r(X_i)$ for $r>0$ and $X_i$ a Poisson point process. ]{{\input{Pix/Fig-Cont} }}%
    \qquad
    \subfloat[\label{Pix_2}Realization of a Poisson--Voronoi tessellation.]{{\input{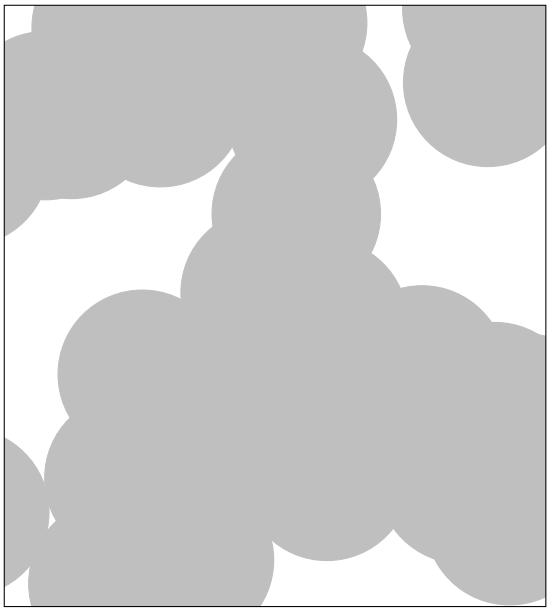}
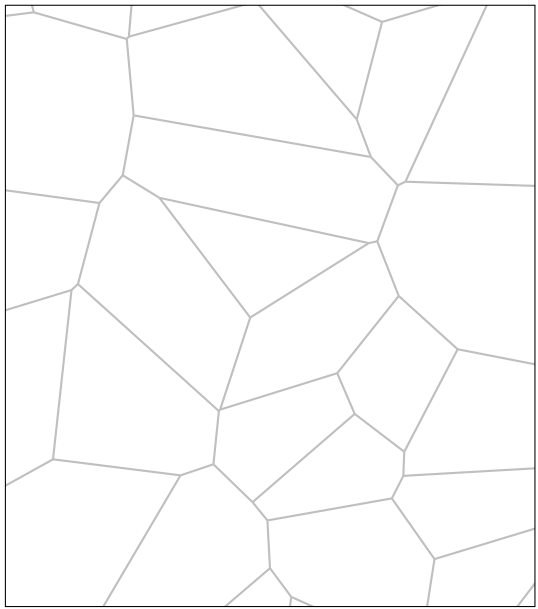}}%
    \caption{Examples of random environments.}%
    \label{fig:example}%
\end{figure}
\begin{example}[Absolutely continuous random environments]\label{Example:Absolutely_Continuous}
    Consider environments of the form $\Sigma(\diff x)= G(x)\diff x$, where $(G(x))_{x\in \R^d}$ is a stationary non-negative random field. For example, 
    \begin{enumerate}[label=(\roman*)]
        \item\label{Example:Random_closed_set} for $\lambda_1,\lambda_2\ge 0$, consider 
        \begin{align*}
            G(x)=\lambda_1\mathds{1}_{\Xi}(x)+\lambda_2\mathds{1}_{\Xi^\compl}(x)
        \end{align*}
        where $\Xi$ is a random closed set, see Figure~\ref{Pix_1} for an illustration. 
        \item consider the \textit{shot-noise field}
        \begin{align*}
            G(x)= \sum_{i\in I}k(x-X_i),
        \end{align*}
        where $(X_i)_{i\in I}$ is a homogeneous Poisson point process and $k$ is integrable.
    \end{enumerate}
    For more information concerning absolutely continuous random environments, see, e.g.,~\cite[Chapter~5]{ChStKe13}.
\end{example}
\begin{example}[Singular random environments]\label{Example:Singular}
    Let $\Sigma(\diff x)= \nu_1(S\cap \diff x)$, where $\nu_1$ is the one-dimensional Hausdorff measure and $S$ is a stationary segment process in the space of line segments. For example, consider
    \begin{enumerate}[label=(\roman*)]
        \item the \textit{Poisson--Voronoi tessellation}, i.e.~for a Poisson point process $(X_i)_{i\in I}$, 
        \begin{align*}
            S=\bigcup_{i,j\in I, \, i\neq j}\{y\in \R^d\,|\, |y-X_i|=|y-X_j|\},
        \end{align*}
        see Figure~\ref{Pix_2} for an illustration. 
        \item the \textit{Poisson--Manhattan grid}~\cite{JaJhVu23}, i.e.~for independent Poisson point processes $(X_i)_{i\in I}$, $(Y_i)_{i\in I}$, 
        \begin{align*}
            S=\bigcup_{i\in I}(\{X_i\}\times \R ) \cup (\R\times \{Y_i\}).
        \end{align*}
    \end{enumerate}
    For more information on singular random environments, see \cite[Chapter 8]{ChStKe13}.
\end{example}
Given an almost sure realization of one of the above random environments $\Sigma$, one can define the Widom--Rowlinson model with environment $\Sigma$ as described in Section~\ref{sec:Setting}. For a visualization of the Widom--Rowlinson model where the environment is given by Poisson--Voronoi tessellation, see Figure~\ref{Pix_3}. 

Applying Corollary~\ref{thm:Random_Environment} to the aforementioned environments yields the following. 
\begin{corollary}\label{Cor:Examples}
    Let $d\ge 2$, $\Sigma$ as in Example~\ref{Example:Absolutely_Continuous}~or~\ref{Example:Singular} and $a>0$. Then, for $z$ large, $\P$-a.s.~$|\mathcal{G}_{z\Sigma, a}|>1$.
\end{corollary}
The corollary above follows by Corollary~\ref{thm:Random_Environment} and known percolation results from the literature \cite{HiJaCa19,JaJhVu23}. See Section~\ref{sec:Random_Environment} for a short review.  
\section{Proofs}\label{sec:Proofs}
The proof of Theorem~\ref{prop:Non-Random_Environment} follows the well-established ideas from \cite{ChChKo95} and is structured as follows. In Section~\ref{sec:Non-Random_Environment}, we show that the finite-volume Widom--Rowlinson models with extremal boundary conditions converge to some measures, which are our candidates for distinct Gibbs measures. Next, in Section~\ref{sec:Random_Cluster_Model}, we identify the difference in densities of these candidates with some expectation w.r.t.~the random cluster model. With the help of FKG domination, as introduced in Section~\ref{Sec:FKG_Domination}, we prove, in Section~\ref{sec:Non_uniqueness}, that this expectation is bounded from below. Ultimately, this proves that the aforementioned candidates from Section~\ref{sec:Non-Random_Environment} are indeed distinct. Meanwhile, Section~\ref{sec:Random_Environment} presents a large class of environments that guarantee percolation, which also includes most of the examples in Corollary~\ref{Cor:Examples}.

\subsection{Existence of infinite-volume Widom--Rowlinson models}\label{sec:Non-Random_Environment}
In order to prove convergence of the finite-volume Widom--Rowlinson model to the infinite one, we need a suitable notion of convergence, namely local convergence. A sequence $(\mu_n)_{n\in \N}$ of probability measures on $\Gamma_\M$ is said to \textit{converge locally to }$\mu$ if 
\begin{align*}
    \lim_{n\to \infty}\int_{\Gamma_\M}f(\pmb{\omega})\, \mu_n(\diff \pmb{\omega})=\int_{\Gamma_\M} f(\pmb{\omega})\,\mu(\diff \pmb{\omega})
\end{align*}
for all $f$ bounded and \textit{local}, i.e.~$f(\pmb{\omega})=f(\pmb{\omega}_\Lambda)$ for some $\Lambda\Subset \M$. Throughout this section, we use the notation $\bar \Lambda=\Lambda\times \{-,+\}$ for $\Lambda\subset\R^d$, $\Lambda_n=[-n,n]^d$ and 
\begin{align*}
        \mu_{\Lambda}^{\pm,\sigma}= \gamma^\sigma_{\bar{\Lambda}}(\cdot\,|\, \partial \Lambda\times \{\pm\})
\end{align*}
as well as $\mu_{n}^{\pm,\sigma}\coloneqq \mu_{\Lambda_n}^{\pm,\sigma}$.
\begin{prop}\label{prop:Exististence_Extreme_Measure}
    There exist measures $\mu^{+,\sigma}$ and $\mu^{-,\sigma }$ such that $\mu_{n}^{\pm,\sigma}$ converges locally to $\mu^{\pm,\sigma }$ along a subsequence $(n_l)$. Furthermore, both $\mu^{+,\sigma}$ and $\mu^{-,\sigma}$ are Gibbs measures of the Widom--Rowlinson model. Also, for all $\Delta\Subset \M$\begin{align}\label{Eq:Convergence_Cardinalities}
        \int_{\Gamma_\M}|\pmb{\omega}_\Delta|\, \mu_{n_l}^{\pm,\sigma}(\diff \pmb{\omega})\xrightarrow{l\to \infty }\int_{\Gamma_\M}|\pmb{\omega}_\Delta|\, \mu^{\pm,\sigma}(\diff \pmb{\omega}).
    \end{align}
\end{prop}
This is essentially~\cite[Theorem~5.7]{Ja18} with the sole difference that we allow boundary conditions $\partial \Lambda_{n}\times \{\pm\}$ that are not in $\Gamma_\M$. To show that the proof remains valid, we highlight its most important steps.
\begin{proof}
    The functions $\rho_m\colon\M^{m}\rightarrow\R$ are called \textit{correlation functions of some probability measure $\mu$} if 
    \begin{align*}
        \int_{\Gamma_\M}F(\pmb{\omega})\, \mu(\diff\pmb{\omega})=\int_{\Gamma_{\M}}f(\overline{ \pmb{x}})\rho_m(\overline{\pmb{x}})\, \sigma^{\otimes m}(\diff \overline{ \pmb{x}}),\quad  \text{for all bounded } f\colon \M^m\rightarrow \R,
    \end{align*}
    where $\sigma^{\otimes m} $ is the $m$-th product measure of $\sigma$ and
    \begin{align*}
        F(\pmb{\omega})=\sum_{\substack{\overline{\pmb{x}}= (\pmb{x}_1,\ldots, \pmb{x}_m)\in \pmb{\omega}^m,\\
        \forall i\neq j: \, \pmb{x}_j\neq \pmb{x}_i}}f(\overline{\pmb{x}}).
    \end{align*}
    By~\cite[Proposition~2.32, Lemma 3.12]{Ja18} $\mu_n^{\pm, \sigma}$ has correlation function
    \begin{align*}
    \rho_m(\pmb{x}_1,\ldots, \pmb{x}_m)=\int_{\Gamma_{\mathbb{M}}} \mathds{1}_{\pmb{V}(\pmb{\zeta}_{\Lambda_n}\cup\partial \Lambda_n \times \{\pm\} )}( \{\pmb{x}_1,\ldots, \pmb{x}_m\} )\, \mu_{n}^{\pm,\sigma}(\diff \pmb{\zeta}),
    \end{align*}
    where 
    \begin{align*}
        \pmb{V}(\pmb{\zeta})=\{\pmb{\omega}\,|\, B_a((\pmb{\omega}\cup \pmb{\zeta})^+)\cap B_a(\pmb{\omega}^-)=\emptyset\text{ and }B_a(\pmb{\omega}^+)\cap B_a((\pmb{\omega}\cup \pmb{\zeta})^-)=\emptyset\}.
    \end{align*}
    In particular, the \textit{Ruelle's bound} as denoted in~\cite[Definition 2.39]{Ja18} is satisfied with $\rho_m\le 1$. By~\cite[Theorem~2.53]{Ja18}, there exist a locally convergent subsequence of $\mu_{n}^{\pm, \sigma}$ with accumulation point $\mu^{\pm, \sigma}$. 

    Furthermore, \cite[Proposition~5.4~(e)]{Ja18} can be extended to $\partial \Lambda_n\times \{\pm\}\notin \Gamma_{\M}$ and, analogously to the proof of~\cite[Theorem 5.6]{Ja18}, by locality of $\gamma^\sigma _{\Delta}(A\,|\, \cdot)$ the accumulations points $\mu^{\pm, \sigma}$ are Gibbs measures of the Widom--Rowlinson model in the environment $\sigma$. Equation \eqref{Eq:Convergence_Cardinalities} follows by~\cite[Theorem~2.51]{Ja18} applied to $\rho_1$ and $f=\mathds{1}_\Delta$.
    \end{proof}
    \noindent As mentioned before, $\mu^{+,\sigma}$ and $\mu^{-,\sigma}$ are our candidates for distinct Gibbs measures as in Theorem~\ref{prop:Non-Random_Environment}. 
    
    Note that a common, alternative approach towards the existence of accumulation points of $\mu_{n}^{\pm,\sigma}$ via level sets of the entropy, compare~\cite[Page 21]{De19}, is not applicable in our setting as it requires stationary environments~$\sigma$. Besides, methods using entropy or Ruelle's bound, a third, common approach towards the existence of $\mu^{\pm,\sigma}$ is based on FKG properties of the two-colored Widom--Rowlinson model \cite{ChChKo95}. However, this procedure yields only a set function $\mu^{\pm,\sigma}$ on the increasing sets and it is not clear to the authors why it can be extended to a measure so easily.
    \subsection{Random cluster model}\label{sec:Random_Cluster_Model}
    \noindent Up to this point, it is not clear if $\mu^{+,\sigma}$ and $\mu^{-,\sigma}$ coincide. The argument goes very closely along the lines of the classical phase transition result when $\sigma$ is the Lebesgue measure. First, we compare the two specifications on some compact set $\Lambda$ and count the average number of particles in some second compact set $\Delta\subset \Lambda$. For this, we denote by $\mathcal{C}_{\partial \Lambda}(\omega)$, for $\omega\in \Gamma_{\R^d}$, the number of connected components in $B_a(\omega\cup\partial \Lambda)$ and define $N_{\Delta,\partial \Lambda}(\omega)$ as the number of particles in $\Delta$ that are in the same component of $B_a(\omega\cup\partial \Lambda)$ as $\partial \Lambda$. Furthermore, we define the \emph{random cluster model} $\nu_{\sigma,\Lambda}$ as a probability measure on $\Gamma_{\R^d}$ defined via 
\begin{align*}
    \nu_{\sigma,\Lambda}(\diff \omega)=\frac{1}{\tilde Z_\Lambda^\sigma}2^{\mathcal{C}_{\partial \Lambda}(\omega )-1}\, \Pi_{\sigma,\Lambda}(\diff\omega ),
\end{align*}
where $\tilde Z_\Lambda^\sigma$ is chosen to normalize. We have the following identity. 
\begin{prop}\label{Proposition_Symmetry_Breaking}
Let $\Delta\subset \Lambda\subset \R^d$ compact, then,
\begin{align*}
\int_{\Gamma_{\mathbb{M}}}|\pmb{\omega}^+_ \Delta|\, \mu_\Lambda^{+,\sigma}(\diff \pmb{\omega})-\int_{\Gamma_{\mathbb{M}}}|\pmb{\omega}^{+}_ \Delta|\, \mu_\Lambda^{-,\sigma}(\diff \pmb{\omega})=\int_{\Gamma_{\R^d}}N_{\Delta,\partial \Lambda}(\omega)\, \nu_{\sigma,\Lambda}(\diff \omega ).
\end{align*}
\end{prop}

\begin{proof}

By symmetry 
\begin{align*}
\Psi&\coloneqq \int_{\Gamma_{\mathbb{M}}}|\pmb{\omega}^+_ \Delta|\, \mu_\Lambda^{+,\sigma}(\diff \pmb{\omega})-\int_{\Gamma_{\mathbb{M}}}|\pmb{\omega}^{+}_ \Delta|\, \mu_\Lambda^{-,\sigma}(\diff \pmb{\omega})\\
&=\int_{\Gamma_{\mathbb{M}}}\left(|\pmb{\omega}^+_ \Delta|-|\pmb{\omega}^-_ \Delta|\right)\, \gamma^{\sigma }_{ \bar \Lambda}(\diff \pmb{\omega}\, |\,\pmb{\omega}^\star ),
\end{align*}
where $\pmb{\omega}^\star=\partial \Lambda\times \{+\}$. Identify some $\overline{x}=(x_1,\ldots, x_n)\in (\R^d)^n$ and $\overline{c}=(c_1,\ldots, c_n)\in \{-,+\}^n$ with some $\pmb{\omega}(\overline{x},\overline{c})\in \Gamma_{\mathbb{M}}$ by writing $\pmb{\omega}(\overline{x},\overline{c})=\{(x_i,c_i)\,|\, i=1,\ldots, n\}$. Then, 
\begin{align*}
\Psi&=\frac{\exp(-2\sigma(\Lambda))}{Z_\Lambda^\sigma(\pmb{\omega}^\star)}\sum_{n\in \N}\sum_{\bar c\in \{-,+\}^n} \frac{1}{n!}\int_{\Lambda^n}\left(\left|\pmb{\omega}(\bar x,\bar c)^{+}_\Delta\right|-\left|\pmb{\omega}(\bar x,\bar c)_\Delta^{-}\right|\right)\mathds{1}_{\pmb{V}}(\pmb{\omega}(\bar x,\bar c)\cup\pmb{\omega} ^\star )\, \sigma^n(\diff \bar x)\\
&\eqcolon \frac{\exp(-2\sigma(\Lambda))}{Z_\Lambda^\sigma(\pmb{\omega}^\star)}\sum_{n\in \N}\frac{1}{n!}\int_{\Lambda ^n}\varphi_n(\bar x)\, \sigma^n(\diff \bar x).
\end{align*}
Fix some $\bar {x}=(x_1,\ldots, x_n)\in (\R^d)^n$ and define $\{\omega_1,\ldots, \omega_{m}\}$ such that 
\begin{enumerate}[label=(\roman*)]
    \item $\bigcup_{j=1}^m\omega_j=\{x_1,\ldots, x_n\}$,
    \item for all $1\le j< m$, the set $\omega_j$ is nonempty and $B_a(\omega_j)$ is a connected component of $B_a(\{x_1,\ldots, x_n\}\cup \partial \Lambda)$, 
    \item $B_a( \omega_m\cup \partial \Lambda)$ is a connected component of $B_a(\{x_1,\ldots, x_n\}\cup \partial \Lambda)$.
\end{enumerate}
Using two times that all points in $\omega_j$ have the same color $+$ or $-$ when restricting to $\pmb{\omega}(\bar x,\bar c)\cup\pmb{\omega}^\star$ that are viable, we arrive at
\begin{align*}
   \varphi_n(\bar x)
&=\sum_{\bar c\in \{-,+\}^n}\sum_{j=1}^m \left(\left|\pmb{\omega}(\bar x,\bar c)_\Delta^{+}\cap \omega_j\right|-\left|\pmb{\omega}(\bar x,\bar c)_\Delta^{-}\cap \omega _j\right|\right)\mathds{1}_{\pmb{V}}(\pmb{\omega}(\bar x,\bar c)\cup\pmb{\omega} ^\star )\\
&=\sum_{\bar c\in \{-,+\}^n}|\omega_m\cap \Delta|\, \mathds{1}_{\pmb{V}}(\pmb{\omega}(\bar x,\bar c)\cup\pmb{\omega} ^\star )+\sum_{( c_1,\ldots, c_{m-1})\in \{-,+\}^{m-1}}\sum_{j=1}^{m-1}c_j|\omega_j\cap \Delta|\\
&=|\omega_m\cap \Delta|\sum_{\bar c\in \{-,+\}^n}\mathds{1}_{\pmb{V}}(\pmb{\omega}(\bar x,\bar c)\cup\pmb{\omega} ^\star )\\
&=|\omega_m\cap \Delta|\,  2^{m-1}\\
&=N_{\Delta,\partial \Lambda}(\{x_1,\ldots, x_n\})\,  2^{\mathcal{C}_{\partial \Lambda}(\{x_1,\ldots, x_n\})-1}.
\end{align*}
Plugging this back into the original equation results in
\begin{align*}
    \Psi=\frac{\exp(-\sigma(\Lambda))}{Z_\Lambda^\sigma(\pmb{\omega}^\star)}\int_{\Gamma_{\R^d}}N_{\Delta,\partial \Lambda}(\omega)2^{\mathcal{C}_{\partial \Lambda}(\omega )-1}\, \Pi_{\Lambda,\sigma}(\diff\omega ).
\end{align*}
With a very similar argument, one can prove that the constant in front normalizes the expression, i.e.,
\begin{align*}
Z_\Lambda^\sigma(\pmb{\omega}^\star)&=\exp(-2\sigma(\Lambda))\sum_{n\in \N}\sum_{\bar c\in \{-,+\}^n} \frac{1}{n!}\int_{\Lambda^n}\mathds{1}_{\pmb{V}}(\pmb{\omega}(\bar x,\bar c)\cup\pmb{\omega} ^\star )\, \sigma^n(\diff \bar x)\\
&=\exp(-2\sigma(\Lambda))\sum_{n\in \N}\frac{1}{n!}\int_{\Lambda^n}2^{\mathcal{C}_{\partial \Lambda}(\{x_1,\ldots, x_n\} )-1}\, \sigma^n(\diff ( x_1,\ldots, x_n))\\
&=\exp(-\sigma(\Lambda))\int_{\Gamma_{\R^d}}2^{\mathcal{C}_{\partial \Lambda}(\omega)-1}\,  \Pi_{\Lambda,\sigma}(\diff\omega )\\
&=\exp(-\sigma(\Lambda)) \tilde Z^{\sigma}_\Lambda,
\end{align*}
as desired. 
\end{proof}
\noindent So, finding a condition under which $\mu^{+,\sigma}$ and $\mu^{-,\sigma}$ do not coincide reduces to estimating the random cluster model. In this context, the concept on FKG domination is useful.

\subsection{FKG domination}\label{Sec:FKG_Domination}
We say that $f\colon\Gamma_\Lambda\rightarrow \R$ measurable is increasing if $ f(\zeta)\le f(\omega)$ for all $\zeta\subset \omega$. Similarly, $A\in \mathfrak{F}_{\Lambda}$ is increasing if $\mathds{1}_A$ is increasing. Finally, we say that a probability measure $Q$ is \textit{stochastically dominated} by $P$, or $Q\preceq P $, if $Q(A)\le P(A)$ for all increasing $A$.
A classical result~\cite{Pr75,GeKu97} connects stochastic domination with the concept of Papangelou intensities. We write $P\in \mathcal{P}_{\sigma,\Lambda}$ if $P$ is a probability measure on $\Gamma_\Lambda$ that is absolutely continuous with respect to $\Pi_{\sigma,\Lambda}$ with density $f_P$ such that $\{f_P=0\}$ is increasing. For $P\in\mathcal{P}_{\sigma,\Lambda} $, we call 
\begin{align*}
    \rho_P(x,\omega)= \begin{cases}
        \frac{f_P(\{x\}\cup \omega)}{f_P(\omega)}& f_P(\omega)\neq 0\\
        0&\text{else}
    \end{cases}
\end{align*}
the \textit{Papangelou intensity} of $P$. With this notion at hand, the main result in~\cite{GeKu97} reads as follows.
\begin{prop}\label{prop:StochDom_Papangelou}
     Let $\Lambda\Subset \R^d$ and $\sigma$ a finite measure on $\Lambda$. Assume that $P,Q\in \mathcal{P}_{\sigma,\Lambda}$ with 
    \begin{align*}
        \rho_Q(x,\omega)\le\rho_P(x,\zeta) \text{ when } \omega\subset \zeta,\, x\notin \zeta\setminus \omega
    \end{align*}
    for all $x\in \Lambda$ and $\omega,\zeta\in \Gamma_{\Lambda}$. Then, $Q$ is stochastically dominated by $P$. 
\end{prop}
Coming back to the setting of the Widom--Rowlinson model, we apply Proposition~\ref{prop:StochDom_Papangelou} to the random cluster model. 
\begin{prop}\label{prop:RandomCluster_Domination}
    There is a constant $\tau>0$ only depending on the dimension $d$ such that, for all $\Lambda\subset\R^d$ compact, $\Pi_{\tau\sigma,\Lambda}$ is stochastically dominated by the random cluster model $\nu_{\sigma, \Lambda}$.
\end{prop}
\begin{proof}
    We will apply Proposition~\ref{prop:StochDom_Papangelou}. The reference measure on $\Lambda$ will be $\tau \sigma$ for some $\tau>0$ that will be specified later. The density $f_{\nu_{\sigma,\Lambda}}$ from $\nu_{\sigma,\Lambda}$ w.r.t.~$\Pi_{\tau\sigma,\Lambda}$ is given by 
\begin{align*}
    f_{\nu_{\sigma,\Lambda}}(\omega )=\tau ^{-|\omega |}C\, 2^{\mathcal{C}_{\partial \Lambda}(\omega)-1}
\end{align*}
for some constant $C>0$. Adding a point to the configuration $\omega$ can reduce the number of connected components only by a fixed number that depends solely on the dimension $d$. Therefore, one can find a $\tau>0$ such that 
\begin{align*}
    \rho_{\nu_{\sigma,\Lambda}}(x,\omega)=\tau^{-1} 2^{\mathcal{C}_{\partial \Lambda}(\omega\cup\{x\})-\mathcal{C}_{\partial \Lambda}(\omega)}\ge 1 =\rho_{\Pi_{\tau\sigma,\Lambda}}(x,\omega).
\end{align*}
The claim follows by Proposition~\ref{prop:StochDom_Papangelou}.
\end{proof}
\subsection{Non-uniqueness of the Widom--Rowlinson model}\label{sec:Non_uniqueness}
The last ingredient towards showing that $\mu^{\pm, \sigma}_\Lambda$ do not coincide is percolation. Recall that $\omega \in \Gamma_{\R^d}$ is $a$-percolating if there is an infinite connected component in $B_a(\omega)$ and introduce the notion of $(a,\Delta)$\textit{-percolating}, which means that an infinite connected component in $B_a(\omega)$ intersects $\Delta$. 
Finally, we are ready to prove Theorem~\ref{prop:Non-Random_Environment}. 
\begin{proof}[Proof of Theorem~\ref{prop:Non-Random_Environment}]
By the prerequisite of \eqref{Cond:Percolation} and continuity of measures from below there is a $ \Delta\Subset \R^d$ such that 
\begin{align*}
    \Pi_{\tau \sigma, \R^d}(\{\omega \,|\,  \omega\text{ is } (a,\Delta)\text{-percolating}\})>0.
\end{align*} 
By Proposition~\ref{Proposition_Symmetry_Breaking},
\begin{align*}
    \Psi(\Lambda)&\coloneqq \int_{\Gamma_{\mathbb{M}}}|\pmb{\omega}^+_ \Delta|\, \mu_\Lambda^{+,\sigma}(\diff \pmb{\omega})-\int_{\Gamma_{\mathbb{M}}}|\pmb{\omega}^{+}_ \Delta|\, \mu_\Lambda^{-,\sigma}(\diff \pmb{\omega})\\
    &\ge \nu_{\sigma,\Lambda}(\{\omega\, |\, \omega\cup \Lambda^\compl\text{ is } (a,\Delta)\text{-percolating}\}).
\end{align*}
Note that $A=\{\omega\, |\, \omega\cup \Lambda^\compl\text{ is } (a,\Delta)\text{-percolating}\}$ is an increasing event. 
By Proposition~\ref{prop:RandomCluster_Domination}, 
\begin{align*}
     \Psi(\Lambda)\ge \Pi_{\tau \sigma, \Lambda}(\{\omega\, |\, \omega\cup \Lambda^\compl\text{ is } (a,\Delta)\text{-percolating}\})\ge \Pi_{\tau \sigma, \R^{d}}(\{\omega\, |\, \omega\text{ is } (a,\Delta)\text{-percolating}\})>0,
\end{align*}
independently of $\Lambda$. Fix the boxes $\Lambda^{(l)}=\Lambda_{n_l}$ as in Proposition~\ref{prop:Exististence_Extreme_Measure}. By Proposition~\ref{prop:Exististence_Extreme_Measure}, we have 
\begin{align*}
    \int_{\Gamma_\mathbb{M}}|\pmb{\omega}^+_\Delta |\, \mu^{+,\sigma}(\diff \pmb{\omega})-\int_{\Gamma_\mathbb{M}}|\pmb{\omega}^+_\Delta |\, \mu^{-,\sigma}(\diff \pmb{\omega})=\lim_{l\to \infty}\Psi(\Lambda^{(l)})>0,
\end{align*}
which finishes the proof. 
\end{proof}
\noindent We close this section by briefly revisiting the question of sufficient conditions for uniqueness. 
\begin{remark}\label{Rem:Uniqueness}
As already mentioned after Corollary~\ref{Cor:Lebesgue_Phase_Transition}, in the case $\sigma=\lambda$, it is proved that for $z$ small enough there is a unique Gibbs measure. To prove a similar result for the inhomogeneous setting of Theorem~\ref{prop:Non-Random_Environment}, one might use a discretization of the percolation event and follow the argument in \cite[Theorem~6.10]{GeHaMa01}. This would require the convergence of the Gibbs measure on $\Lambda_n\Subset\R^d$ with empty boundary condition, i.e., a stronger version of Proposition~\ref{prop:Exististence_Extreme_Measure}. To get this convergence along the whole sequence not just a subsequence, one needs an FKG property of the random cluster model $\nu_\Lambda$ with empty boundary conditions for varying $\Lambda$ together with a continuous version of the Edwards--Sokal coupling, i.e., a variant of Proposition~\ref{Proposition_Symmetry_Breaking}. 
\end{remark}

\subsection{Widom--Rowlinson models in random environments}\label{sec:Random_Environment}
The main goal of this subsection is to formulate a class of random environments, namely stabilizing, asymptotically essentially connected random environments, that contains the examples of Section \ref{sec:Examples}, except from the Poisson-Manhattan grid. In doing so, we follow the notation and statements of \cite{HiJaCa19}.
A random measure $\Sigma$ is called \textit{stabilizing} with \textit{stabilization radii} $(R_x)_{x \in\R^d}$ if 
\begin{enumerate}[label=(\roman*)]
    \item \label{Def_Stability_Stationary}$(\Sigma, R)$ is defined on $(\Omega,\mathcal{F},\P)$ and jointly stationary, i.e.~for all $y\in \R^d$ the distributions of $(y+\Sigma, (R_{x+y})_{x\in \R^d})$ coincide, where $(x+\Sigma)(A)=\Sigma(x+A)$ and $x+A=\{x+y\, |\,y\in A\}$,
    \item \label{Def_Stability_Radii}$\lim_{n\to \infty}\P(\sup_{y\in Q_n(0)\cap \Q^d}R_y<n)=1$, where $Q_n(x)=x+[-n/2,n/2]^d$, and 
    \item \label{Def_Stability_Independence}for $n\ge 1$, for all $x_1,\ldots , x_m\in \R^d$ with $|x_i-x_j|>3n$ when  $i\neq j$ and all bounded measurable functions $f\colon\mathcal{M}\rightarrow \R$, 
    \begin{align*}
    \left(f(\Sigma_{Q_n(x_i)})\mathds{1}\{\sup_{y\in Q_n(x_i)\cap \Q^d}R_y<n\}\right)_{i=1,\ldots, m}
    \end{align*}
    is a family of independent random variables.
\end{enumerate}
\begin{remark}\label{rem:Stability_ergodicity}
    Stabilization can be seen as a weak version of assuming that $\Sigma_{A}$ and $\Sigma_B$ are independent for $A$ and $B$ disjoint. Similar to independence, stabilization implies ergodicity of $\Sigma$, i.e.~$\P(\Sigma\in A)\in \{0,1\}$ for all $A\in \mathscr{F}$ with $\{x+\Sigma \in A\}=\{\Sigma \in A\}$ for all $x\in \R^d$. 
\end{remark}
\begin{proof}
Fix $A,B\in \mathscr{F}_\Lambda$ for some $\Lambda\Subset\R^d$ and  $x\in \R^d$ with $|x|=1$. Then, for $n$ so large that $\Lambda\subset B_n(0)$
\begin{align*}
    |\P(\Sigma \in A,\, (4nx+\Sigma)\in B)-\P(\Sigma \in A)\P(\Sigma \in B)|&\le 4\P(\sup_{y\in Q_n(0)\cap \Q^d}R_y\ge n)\xrightarrow{n\to \infty }0.
\end{align*}
Since $\mathscr{F}_\Lambda$ generates $\mathscr{F}$, this convergence holds also for all $A,B\in \mathscr{F}$. Ergodicity follows by choosing $A=B$ with $\{x+\Sigma \in A\}=\{\Sigma \in A\}$ for all $x\in \R^d$, since then $\P(\Sigma\in A)=\P(\Sigma\in A)^2$.
\end{proof}
Apart from some weak independence assumption, we also need some connectedness property to guarantee percolation. A stabilizing measure $\Sigma$ with stabilization radii $(R_x)$ is \textit{assymptotically essentially connected} if, whenever $\sup_{y\in Q_{2n}}R_y<n/2$
\begin{enumerate}[label=(\roman*)]
    \item $\supp(\Sigma_{Q_n(0)})\neq \emptyset$, where 
    \begin{align*}
        \supp(\sigma)=\{x\in \R^d\,|\, \forall \varepsilon>0:\, \sigma(Q_\varepsilon(x))>0\},
    \end{align*}
    \item $\supp(\Sigma_{Q_n(0)})$ is contained in a connected component of $\supp(\Sigma_{Q_{2n}(0)})$.
\end{enumerate}
Indeed, asymptotically essentially connected, stabilizing environments guarantee a supercritical regime of percolation for the corresponding Cox point process. 
\begin{prop}\label{Prop:Stabilizing_Percolation}   Fix an asymptotically essentially connected, stabilizing random measure $\Sigma\colon(\Omega, \mathcal{F},\P)\rightarrow(\mathcal{M},\mathscr{F})$ and the corresponding distribution $\Q_z(A)=\E_\P[\Pi_{z\Sigma,\R^d}(A)]$ of the Cox point process. Then, for all $a>0$ and for sufficiently large $z$,
\begin{align*}
    \Q_z(\{\omega\, |\, \omega \text{ is }a\text{-percolating}\})=1.
\end{align*}
\end{prop}
\begin{proof}
    By~\cite[Proof of Theorem 2.6]{HiJaCa19}, we get that the event of percolation has positive probability for $z$ sufficiently large. By \cite[Proposition 12.3.VII]{DaVe07} ergodicity of the environments implies ergodicity of the Cox point process. So, Remark~\ref{rem:Stability_ergodicity} implies that the probability of percolation is $1$.  
\end{proof}
Finally, we can prove Corollary~\ref{Cor:Examples}.
\begin{proof}[Proof of Corollary~\ref{Cor:Examples}]
By Proposition~\ref{Prop:Stabilizing_Percolation}, Corollary~\ref{thm:Random_Environment} and the examples for stabilizing, asymptotically essentially connected random measures given in \cite{HiJaCa19}, we only need to prove the claim for the non-stabilizing Poisson-Manhattan grid. Here, \cite{JaJhVu23} proves percolation for $z$ large enough. Ergodicity of the Poisson-Manhattan grid follows analogously to Remark~\ref{rem:Stability_ergodicity} since for $A,B\in \mathscr{F}_\Lambda$ where $\Lambda\Subset\R^d$ and for $x=(1,1)$ it holds for sufficiently large $n$
\begin{align*}
    \P(\Sigma \in A,nx+\Sigma\in B )=\P(\Sigma \in A)\P(\Sigma\in B ).
\end{align*}
So, again by \cite[Proposition 12.3.VII]{DaVe07} also the corresponding Cox point process is ergodic and Corollary~\ref{thm:Random_Environment} proves the claim.
\end{proof}
\section*{Acknowledgement} 
BJ received support by the Leibniz Association within the Leibniz Junior Research Group on \textit{Probabilistic Methods for Dynamic Communication Networks} as part of the Leibniz Competition (grant no.\ J105/2020). The research of BJ and DK is funded by Deutsche Forschungsgemeinschaft (DFG) through the SPP2265 within the Project P27.

\section*{References}
\renewcommand*{\bibfont}{\footnotesize}
\printbibliography[heading = none]


\end{document}

%% file: Pix/Fig-Cont.tex
\begin{tikzpicture}[scale=1.6] 
 \begin{scope} 
\clip(5.5,0) rectangle (10,5.0);
\draw[lightgray , fill] (8.970665537195522,1.9680428807852413) circle (7.0mm);
\draw[lightgray , fill] (0.5122982888761973,1.3481198891426782) circle (7.0mm);
\draw[lightgray , fill] (1.170285061164541,3.6472476070059114) circle (7.0mm);
\draw[lightgray , fill] (3.3297925892953186,4.784198741605015) circle (7.0mm);
\draw[lightgray , fill] (7.635756850683907,1.9886053041744545) circle (7.0mm);
\draw[lightgray , fill] (6.4027138274224775,0.19137010762836693) circle (7.0mm);
\draw[lightgray , fill] (8.060383716924669,4.048926519851288) circle (7.0mm);
\draw[lightgray , fill] (0.036038843698021905,4.652334156763529) circle (7.0mm);
\draw[lightgray , fill] (7.520743603087733,4.649961781170127) circle (7.0mm);
\draw[lightgray , fill] (7.039006552432669,0.38721346767512266) circle (7.0mm);
\draw[lightgray , fill] (6.797012767632207,4.189796190784996) circle (7.0mm);
\draw[lightgray , fill] (3.3135133487201704,4.752617907753742) circle (7.0mm);
\draw[lightgray , fill] (4.18605800828462,2.0466254825117147) circle (7.0mm);
\draw[lightgray , fill] (5.841526855481661,4.082445163925856) circle (7.0mm);
\draw[lightgray , fill] (9.698991937653856,0.7137022569955509) circle (7.0mm);
\draw[lightgray , fill] (6.873390029701252,4.710024785767715) circle (7.0mm);
\draw[lightgray , fill] (7.81309121614523,4.852748607455015) circle (7.0mm);
\draw[lightgray , fill] (3.235986926244884,1.374493749586545) circle (7.0mm);
\draw[lightgray , fill] (7.926498433163336,3.264839213089724) circle (7.0mm);
\draw[lightgray , fill] (1.056547172956468,2.2004639720674044) circle (7.0mm);
\draw[lightgray , fill] (6.646805714538835,1.9318416057479748) circle (7.0mm);
\draw[lightgray , fill] (3.0234505500505082,3.7130870026938156) circle (7.0mm);
\draw[lightgray , fill] (4.0704283646483965,3.574527044920053) circle (7.0mm);
\draw[lightgray , fill] (7.667660439946182,2.6082418540948673) circle (7.0mm);
\draw[lightgray , fill] (1.7293286856056744,4.726983464361824) circle (7.0mm);
\draw[lightgray , fill] (0.9107030481693645,0.04110485069153236) circle (7.0mm);
\draw[lightgray , fill] (5.1718084806061775,0.7844930740270828) circle (7.0mm);
\draw[lightgray , fill] (4.659923785137677,0.2481596347386128) circle (7.0mm);
\draw[lightgray , fill] (8.104024253615442,1.6282219412576682) circle (7.0mm);
\draw[lightgray , fill] (6.430092412146702,4.8129306002609065) circle (7.0mm);
\draw[lightgray , fill] (8.180907746278137,1.082467156325485) circle (7.0mm);
\draw[lightgray , fill] (8.152680461520164,2.358498090589134) circle (7.0mm);
\draw[lightgray , fill] (0.33716528192896833,0.4412111139287256) circle (7.0mm);
\draw[lightgray , fill] (1.696767205417623,2.3746002523196843) circle (7.0mm);
\draw[lightgray , fill] (9.575013761314747,1.6629329014354577) circle (7.0mm);
\draw[lightgray , fill] (1.8240508997401006,4.882742640618684) circle (7.0mm);
\draw[lightgray , fill] (6.535838399191104,1.068767515049035) circle (7.0mm);
\draw[lightgray , fill] (9.805532446458997,1.5743754492694007) circle (7.0mm);
\draw[lightgray , fill] (7.248746160880389,4.8772841112873335) circle (7.0mm);
\draw[lightgray , fill] (3.2352382750465924,2.4828028906114463) circle (7.0mm);
\draw[lightgray , fill] (6.05864036749422,4.09210741763089) circle (7.0mm);
\draw[lightgray , fill] (3.558740515137556,3.8010466830741096) circle (7.0mm);
\draw[lightgray , fill] (4.24258848416005,1.3396220722112795) circle (7.0mm);
\draw[lightgray , fill] (1.8289522374541167,1.80288625100391) circle (7.0mm);
\draw[lightgray , fill] (3.9752797959445685,3.582825849548371) circle (7.0mm);
\draw[lightgray , fill] (3.8458477364191634,2.8392708438707026) circle (7.0mm);
\draw[lightgray , fill] (0.05550731753927063,3.785848799023701) circle (7.0mm);
\draw[lightgray , fill] (9.715941716517978,4.594801402536982) circle (7.0mm);
\draw[lightgray , fill] (8.196628871317987,1.3429867250695138) circle (7.0mm);
\draw[lightgray , fill] (0.8927304178716666,4.450777612426328) circle (7.0mm);
\draw[lightgray , fill] (6.682034440850952,0.23853801283035714) circle (7.0mm);
\draw[lightgray , fill] (9.303123561281875,1.0581403301885195) circle (7.0mm);
\draw[lightgray , fill] (2.873804029522674,0.3700251673733146) circle (7.0mm);
\draw[lightgray , fill] (3.8731427244966166,3.07965773132275) circle (7.0mm);
\draw[lightgray , fill] (7.091923037868739,4.615076080688394) circle (7.0mm);
\draw[lightgray , fill] (4.653250833494217,1.3043739354285633) circle (7.0mm);
\draw[lightgray , fill] (9.516890002081146,4.358530777751298) circle (7.0mm);
\draw[lightgray , fill] (5.174198517516223,3.699378301510822) circle (7.0mm);
\draw[lightgray , fill] (3.763854635688144,4.250611943738388) circle (7.0mm);
\draw[lightgray , fill] (2.654414946159699,1.3661189945169294) circle (7.0mm);
\draw[lightgray , fill] (7.617509371612277,1.4388025009660266) circle (7.0mm);
\draw[lightgray , fill] (2.88054165619674,4.617649544308659) circle (7.0mm);
\draw[lightgray , fill] (9.510025053283064,4.97072347789058) circle (7.0mm);
\draw[lightgray , fill] (1.9342067237671123,0.7048334668166301) circle (7.0mm);
\draw[lightgray , fill] (1.1297568740355268,1.6335820759647368) circle (7.0mm);
 \end{scope} 
 \draw (5.5,0) rectangle (10,5.0); 
 \end{tikzpicture}

%% file: refs/refs.bib
@incollection{GeHaMa01,
title = {The random geometry of equilibrium phases},
editor = {C. Domb and J.L. Lebowitz},
series = {Phase Transitions and Critical Phenomena},
publisher = {Academic Press},
volume = {18},
pages = {1-142},
year = {2001},
issn = {1062-7901},
author = {Hans-Otto Georgii and Olle Häggström and Christian Maes}
}

@article{Ru71,
  title = {Existence of a Phase Transition in a Continuous Classical System},
  author = {Ruelle, David},
  journal = {Physical Review Letters},
  volume = {27},
  pages = {1040--1041},
  numpages = {0},
  year = {1971},
publisher = {American Physical Society}
}

@article{Ja18,
  title={Gibbsian point processes},
  author={Jansen, Sabine},
  journal={available at author’s website},
  year={2018}
}

@article{JaJhVu23,
author = {Benedikt Jahnel and Sanjoy Kumar Jhawar and Anh Duc Vu},
title = {{Continuum percolation in a nonstabilizing environment}},
volume = {28},
journal = {Electronic Journal of Probability},
publisher = {Institute of Mathematical Statistics and Bernoulli Society},
pages = {1 -- 38},
year = {2023},
doi = {10.1214/23-EJP1029}
}

@book{BoRi06, place={Cambridge}, title={Percolation}, publisher={Cambridge University Press}, author={Bollobás, Bela and Riordan, Oliver}, year={2006}}

@book{Me96, place={Cambridge}, series={Cambridge Tracts in Mathematics}, title={Continuum Percolation}, publisher={Cambridge University Press}, author={Meester, Ronald and Roy, Rahul}, year={1996}, collection={Cambridge Tracts in Mathematics}}

@article{Pr75,
  title={Spatial birth and death processes},
  author={Chris Preston},
  journal={Advances in Applied Probability},
  year={1975},
  volume={7},
  pages={465 - 466}
}

@article{DeHo19,
  title={Phase transition for continuum Widom-Rowlinson model with random radii},
  author={Dereudre, David and Houdebert, Pierre},
  journal={Journal of Statistical Physics},
  volume={174},
  pages={56--76},
  year={2019},
  publisher={Springer}
}

@article{DeHo21,
author = {David Dereudre and Pierre Houdebert},
title = {{Sharp phase transition for the continuum Widom-Rowlinson model}},
volume = {57},
journal = {Annales de l'Institut Henri Poincar\'e, Probabilit\'es et Statistiques},
number = {1},
publisher = {Institut Henri Poincar\'e},
pages = {387 -- 407},
year = {2021}
}

@article{DiLiXi24,
  title={Long range order for three-dimensional random field Ising model throughout the entire low temperature regime},
  author={Ding, Jian and Liu, Yu and Xia, Aoteng},
  journal={Inventiones Mathematicae},
  volume={238},
  number={1},
  pages={247--281},
  year={2024},
  publisher={Springer}
}

@article{Ch24,
  title={Features of a spin glass in the random field Ising model},
  author={Chatterjee, Sourav},
  journal={Communications in Mathematical Physics},
  volume={405},
  number={493},
  year={2024},
  publisher={Springer}
}

@article{BrKu88,
  title={Phase transition in the 3d random field Ising model},
  author={Bricmont, Jean and Kupiainen, Antti},
  journal={Communications in Mathematical Physics},
  volume={116},
  pages={539--572},
  year={1988},
  publisher={Springer}
}

@article{FyMaPiSo18,
  title={Review of recent developments in the random-field Ising model},
  author={Fytas, Nikolaos G and Mart{\'\i}n-Mayor, V{\'\i}ctor and Picco, Marco and Sourlas, Nicolas},
  journal={Journal of Statistical Physics},
  volume={172},
  pages={665--672},
  year={2018},
  publisher={Springer}
}

@inproceedings{Ku19,
  title={Gibbs-non Gibbs transitions in different geometries: The Widom-Rowlinson model under stochastic spin-flip dynamics},
  author={K{\"u}lske, Christof},
  booktitle={Statistical Mechanics of Classical and Disordered Systems: Luminy, France, August 2018},
  pages={3--19},
  year={2019},
  organization={Springer}
}

@article{HiTa04,
author = {Yasunari Higuchi and Masato Takei},
title = {{Some results on the phase structure of the two-dimensional Widom-Rowlinson model}},
volume = {41},
journal = {Osaka Journal of Mathematics},
number = {2},
publisher = {Osaka University and Osaka Metropolitan University, Departments of Mathematics},
pages = {237 -- 255},
year = {2004},
}

@article{HiJaMu22,
  title={Sharp phase transition for Cox percolation},
  author={Hirsch, Christian and Jahnel, Benedikt and Muirhead, Stephen},
  journal={Electronic Communications in Probability},
  volume={27},
  pages={1--13},
  year={2022},
  publisher={The Institute of Mathematical Statistics and the Bernoulli Society}
}

@article{JaToCa22,
  title={Phase transitions for the Boolean model of continuum percolation for Cox point processes},
  author={Jahnel, Benedikt and T{\'o}bi{\'a}s, Andr{\'a}s and Cali, Elie},
  journal={Brazilian Journal of Probability and Statistics},
  volume={36},
  number={1},
  pages={20--44},
  year={2022},
  publisher={Brazilian Statistical Association}
}

@article{WiRo70,
  title={New model for the study of liquid-vapor phase transitions},
  author={Widom, Benjamin and Rowlinson, John S},
  journal={The Journal of Chemical Physics},
  volume={52},
  number={4},
  pages={1670--1684},
  year={1970},
  publisher={AIP Publishing}
}

@article{DeHo15,
author = {David Dereudre and Pierre Houdebert},
title = {{Infinite volume continuum random cluster model}},
volume = {20},
journal = {Electronic Journal of Probability},
number = {none},
publisher = {Institute of Mathematical Statistics and Bernoulli Society},
pages = {1 -- 24},
year = {2015}
}

@article{GiLeMa95,
  title={Agreement percolation and phase coexistence in some Gibbs systems},
  author={Giacomin, G and Lebowitz, JL and Maes, Ch},
  journal={Journal of Statistical Physics},
  volume={80},
  pages={1379--1403},
  year={1995},
  publisher={Springer}
}

@book{DaVe07,
  title={An Introduction to the Theory of Point Processes: Volume II: General Theory and Structure},
  author={Daley, D.J. and Vere-Jones, D.},
  series={Probability and Its Applications},
  year={2007},
  publisher={Springer New York}
}

@article{HiJaCa19,
title = {Continuum percolation for Cox point processes},
journal = {Stochastic Processes and their Applications},
volume = {129},
number = {10},
pages = {3941-3966},
year = {2019},
issn = {0304-4149},
author = {Christian Hirsch and Benedikt Jahnel and Elie Cali}
}

@article{GeKu97,
 ISSN = {00219002},
 author = {Hans-Otto Georgii and Torsten Küneth},
 journal = {Journal of Applied Probability},
 number = {4},
 pages = {868--881},
 publisher = {Applied Probability Trust},
 title = {Stochastic comparison of point random fields},
 volume = {34},
 year = {1997}
}

@article{ChChKo95,
author = {J. T. Chayes and L. Chayes and R. Koteck{\'y}},
title = {{The analysis of the Widom-Rowlinson model by stochastic geometric methods}},
volume = {172},
journal = {Communications in Mathematical Physics},
number = {3},
publisher = {Springer},
pages = {551 -- 569},
year = {1995},
}

@article{GeHa95,
  title={Phase transition in continuum Potts models},
  author={Georgii, Hans-Otto and H{\"a}ggstr{\"o}m, Olle},
  journal={Communications in Mathematical Physics},
  volume={181},
  number={2},
  pages={507--528},
  year={1996},
  publisher={Springer}
}

@book{GuPe17, place={Cambridge}, series={Institute of Mathematical Statistics Textbooks}, title={Lectures on the Poisson Process}, publisher={Cambridge University Press}, author={Last, Günter and Penrose, Mathew}, year={2017}, collection={Institute of Mathematical Statistics Textbooks}}

@book{ChStKe13,
  title={Stochastic Geometry and Its Applications},
  author={Chiu, S.N. and Stoyan, D. and Kendall, W.S. and Mecke, J.},
  series={Wiley Series in Probability and Statistics},
  year={2013},
  publisher={Wiley}
}

@Inbook{De19,
author="Dereudre, David",
title="Introduction to the theory of Gibbs point processes",
bookTitle="Stochastic Geometry: Modern Research Frontiers",
year="2019",
publisher="Springer International Publishing",
address="Cham",
pages="181--229"
}

@book{JaKo20,
  title={Probabilistic Methods in Telecommunications},
  author={Jahnel, Benedikt and K{\"o}nig, Wolfgang},
    series={Compact Textbooks in Mathematics},
  year={2020},
  publisher={Springer}
}

@book{Ge11,
  title={Gibbs Measures and Phase Transitions},
  author={Georgii, Hans-Otto},
  volume={9},
    series={De Gruyter Studies in Mathematics},
  year={2011},
  publisher={Walter de Gruyter}
}
